\documentclass[final]{elsarticle} 

\usepackage{amssymb,latexsym}
\usepackage{amsmath}
\usepackage{graphicx}
\usepackage{lscape}
\biboptions{sort&compress}

\journal{: \: Indagationes Mathematicae}

\begin{document}

\newtheorem{lema}{Lemma}
\newtheorem{teo}{Theorem}
\newproof{prova}{Proof}

\begin{frontmatter}

\title{New definite integrals and a two-term dilogarithm identity}

\author{F. M. S. Lima}

\address{Institute of Physics, University of Brasilia, P.O. Box 04455, 70919-970, Brasilia, DF, Brazil}


\ead{fabio@fis.unb.br}

\date{\today}

\begin{abstract}
Among the several proofs known for $\sum_{n=1}^\infty{\! 1/n^2} = {\pi^2/6}$, the one by Beukers, Calabi, and Kolk involves the evaluation of $\,\int_0^1\!\!{\int_0^1{1/(1-x^2 y^2) \, dx} \, dy}$. It starts by showing that this double integral is equivalent to $\frac34 \sum_{n=1}^\infty{\! 1/n^2}$, and then a non-trivial \emph{trigonometric} change of variables is applied which transforms that integral into $\,{\int\!\!\int}_T \: 1 \; du \, dv$, where $T$ is a triangular domain whose area is simply ${\pi^2/8}$.  Here in this note, I introduce a hyperbolic version of this change of variables and, by applying it to the above integral, I find exact closed-form expressions for $\int_0^\infty{\left[\sinh^{-1}{(\cosh{u})}-u\right] d u}$, $\,\int_{\alpha}^\infty{\left[u-\cosh^{-1}{(\sinh{u})}\right] d u}$, and $\,\int_{\,\alpha/2}^\infty{\ln{(\tanh{u})} \: d u}$, where $\alpha = \sinh^{-1}(1)$.  From the latter integral, I also derive a two-term dilogarithm identity.
\newline
\end{abstract}


\end{frontmatter}

\section{Introduction}
The Riemann zeta function $\zeta{(s)}$ is defined, for complex values of $s$, $\Re{(s)}>1$, by $\zeta{(s)} := \sum_{n=1}^\infty{\,1/n^s}$.  For integer values of $s$, $s>1$, the first zeta value is $\zeta{(2)} = \sum_{n=1}^\infty{1/n^2}$. Euler was the first to estimate the numerical value of this series to more than $5$ decimal places (1735), as well as to determine an exact closed-form expression for it (1740), as given by~\cite{Euler}:

\begin{lema} [Euler's result for $\zeta{(2)}$]
\label{lem:z2}
\quad $\sum_{n=1}^\infty{\! 1/n^2} = {\pi^2/6}$.
\end{lema}

From Euler's time onward, this result has been proved in several forms, from elementary to complex ones~\cite{Kalman,Chapman}. Among these proofs, at least two involve the evaluation of a double integral over the unit square and these are the ones in which we are interested here. Let us briefly describe these two proofs for a better exposition of the hyperbolic approach and its consequences.

\begin{prova}[Apostol's proof]
The short proof by Apostol (1983) can be found in Ref.~\cite{Apostol}. He starts by expanding the integrand of Beukers' integral $I := \int_0^1\!\!\!{\int_0^1{\!\frac{1}{1-x\,y}~dx} \, dy}$ as a geometric series, in order to show that it is equivalent to
\begin{equation*}
\int_0^1\!\!\!{\int_0^1{\left(1 +x y +x^2 y^2 +\ldots\right) dx} \, dy} = \int_0^1{\left(1+\frac{y}{2}+\frac{y^2}{3}+\ldots\right) dy} 
= \sum_{n=1}^\infty{\frac{1}{n^2}} \, ,
\end{equation*}
the interchange of integrals and sums being justifiable by the fact that each integrand is nonnegative and the sums converge absolutely (see the only corollary of Theorem~7.16 in Ref.~\cite{bk:Rudin}).\footnote{Note that ${\,1/(1-x y)}$ is a positive real number at every point in the unit square, the only exception being the point $(1,1)$, where it is undefined. Of course, usual methods for treating improper integrals apply to this case. \medskip} By applying a simple linear change of variables corresponding to the rotation of the coordinate axes through the angle $\pi/4$~radians, as given by\footnote{For this change of variables, it is easy to see that $\left|J\right| = 1$, $J$ being the Jacobi determinant $\frac{d(x,y)}{d(u,v)} = \left|
  \begin{array}{cc}
    {\partial x/\partial u} & {\partial x/\partial v} \\
    {\partial y/\partial u} & {\partial y/\partial v} \\
  \end{array} \right|$.}
\begin{equation*}
x = \frac{u-v}{\sqrt{2}} \quad \text{and} \quad y = \frac{u+v}{\sqrt{2}} \, ,
\end{equation*}
directly to the Beukers' integral, one finds~\cite{Apostol,bk:TheBook}:
\begin{equation*}
I = 4 \int_0^{\sqrt{2}/2}\:{\frac{\arctan\!{\left( \frac{u}{\sqrt{2-u^2}}\right)}}{\sqrt{2-u^2}}\, du} 
+ 4 \int_{\,\sqrt{2}/2}^{\sqrt{2}}{\frac{\arctan\!{\left( \frac{\sqrt{2}-u}{\sqrt{2-u^2}}\right)}}{\sqrt{2-u^2}}\, du} \, .
\end{equation*}
Now, the substitution $\,u = \sqrt{2} \: \sin{\theta}$, followed by the use of the trigonometric identity $\,\arctan{\left(\sec{\theta}-\tan{\theta}\right)} = {\,\pi/4}-{\,\theta/2}\,$ on the second integral, yields $I = {\,\pi^2/18} + {\,\pi^2/9}$. Therefore, $\sum_{n=1}^\infty{1/n^2} = {\,\pi^2/18} + {\,\pi^2/9} = {\,\pi^2/6}$.  \qquad \qquad \qquad $\Box$ \end{prova}

The other proof, in which we are more interested here, was given by Beukers, Calabi, and Kolk (1993)~\cite{Calabi}.

\begin{prova}[BCK proof]
Similarly to Apostol's proof, Beukers, Calabi, and Kolk (BCK) start by showing that
\begin{eqnarray}
K &:=& \int_0^1\!\!\!{\int_0^1{\frac{1}{1-x^2 y^2} \, dx} \, dy} = \int_0^1\!\!\!{\int_0^1{\sum_{n=0}^\infty{x^{2n} y^{2n}} \, dx} \, dy} \nonumber \\
&=& \sum_{n=0}^\infty{\int_0^1\!\!\!{\int_0^1{x^{2n} y^{2n}} \, dx} \, dy} = \sum_{n=0}^\infty{\frac{1}{(2n+1)^2}} \, .
\label{eq:intK}
\end{eqnarray}
Again, as each integrand is nonnegative the absolute convergence allows the interchange of the sum and the integral. Since $\,\sum_{n=1}^\infty{{\,1/(2n)^2}} = \frac14 \, \sum_{n=1}^\infty{{\,1/n^2}}$, then $\:\sum_{n=0}^\infty{\frac{1}{(2n+1)^2}} = \frac34 \, \sum_{n=1}^\infty{\frac{1}{n^2}}$.
They then evaluate the unit square integral that defines $K$ by applying the following non-trivial, trigonometric change of variables~\cite{Calabi}:
\begin{equation*}
x = \frac{\sin{u}}{\cos{v}} \; , \quad y = \frac{\sin{v}}{\cos{u}} \: .
\end{equation*}
For this change of variables, one easily finds that $|J|= 1-\tan^2{u}\,\tan^2{v}=1-x^2 y^2$, which reduces the integral for $K$ to $\, {\int \! \! \! \int}_{\!\! T} \: 1 \; du \, dv \,$
where $T$ is the domain in the $u v$-plane defined by $\{(u,v): \, u \ge 0$, $v \ge 0$, $u+v \le {\,\pi/2} \}$~\cite{Kalman,bk:TheBook,Calabi}. This double integral exactly evaluates to ${\,\pi^2/8\,}$, i.e. the area of the triangle corresponding to $\,T$, which implies Lemma~\ref{lem:z2}. \qquad \qquad \qquad \qquad \qquad \qquad \qquad $\Box$
\end{prova}

Here in this work, I introduce a hyperbolic version of BCK change of variables for double integrals over the unit square. By applying this new change of variables to $\int_0^1\!\!{\int_0^1{1/(1-x^2 y^2) \, dx} \, dy}$, I derive exact closed-form expressions for some definite integrals.  From one of these integrals, I deduce a closed-form expression for $\,\mathrm{Li}_2{(\sqrt{2}-1)} + \mathrm{Li}_2{\left(1-{\,\sqrt{2}/2}\right)}$.

\section{A hyperbolic version of BCK change of variables}
\begin{figure}[!h]
\begin{center}
\scalebox{0.16}{\includegraphics{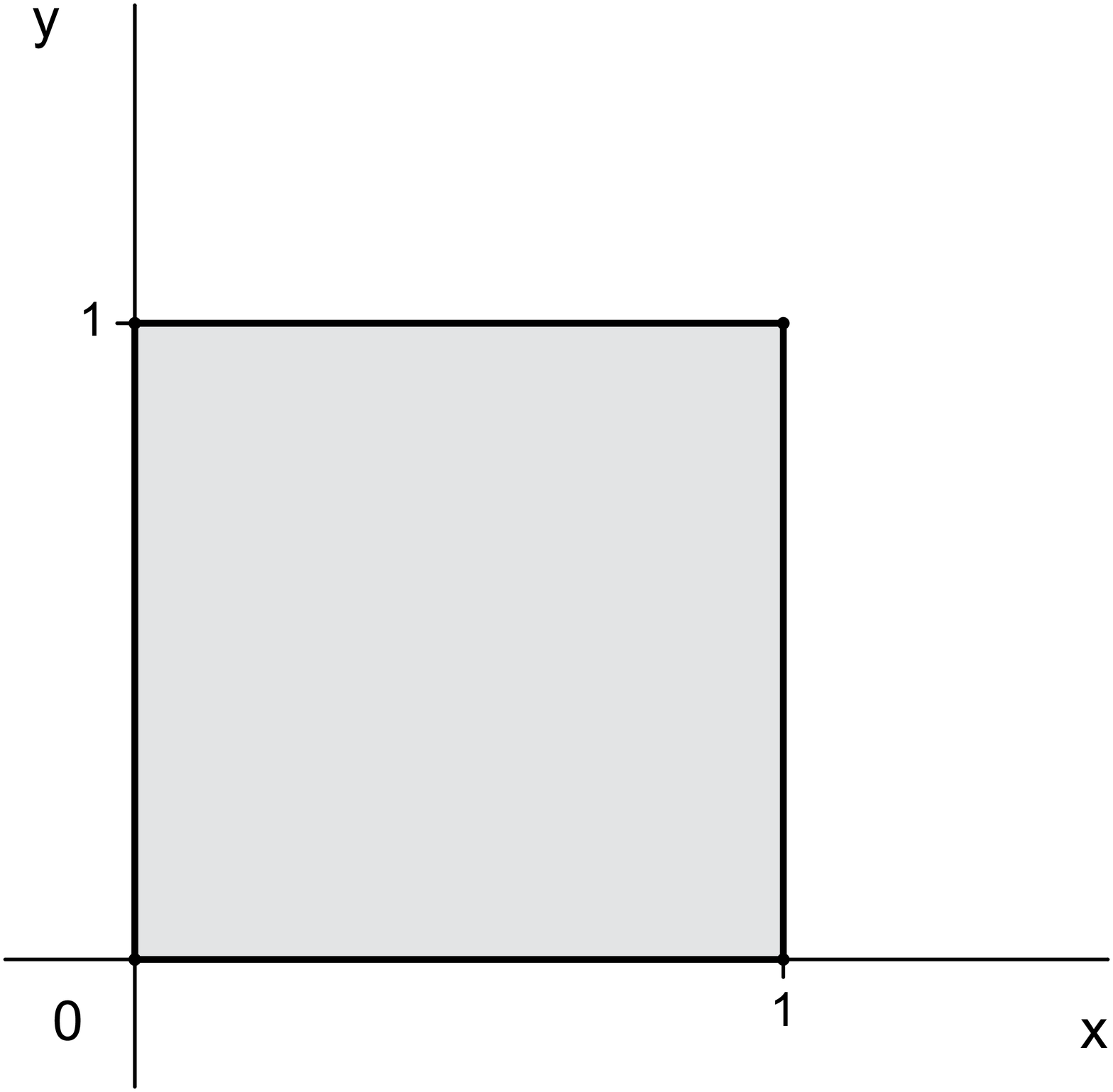}}
\scalebox{0.27}{\includegraphics{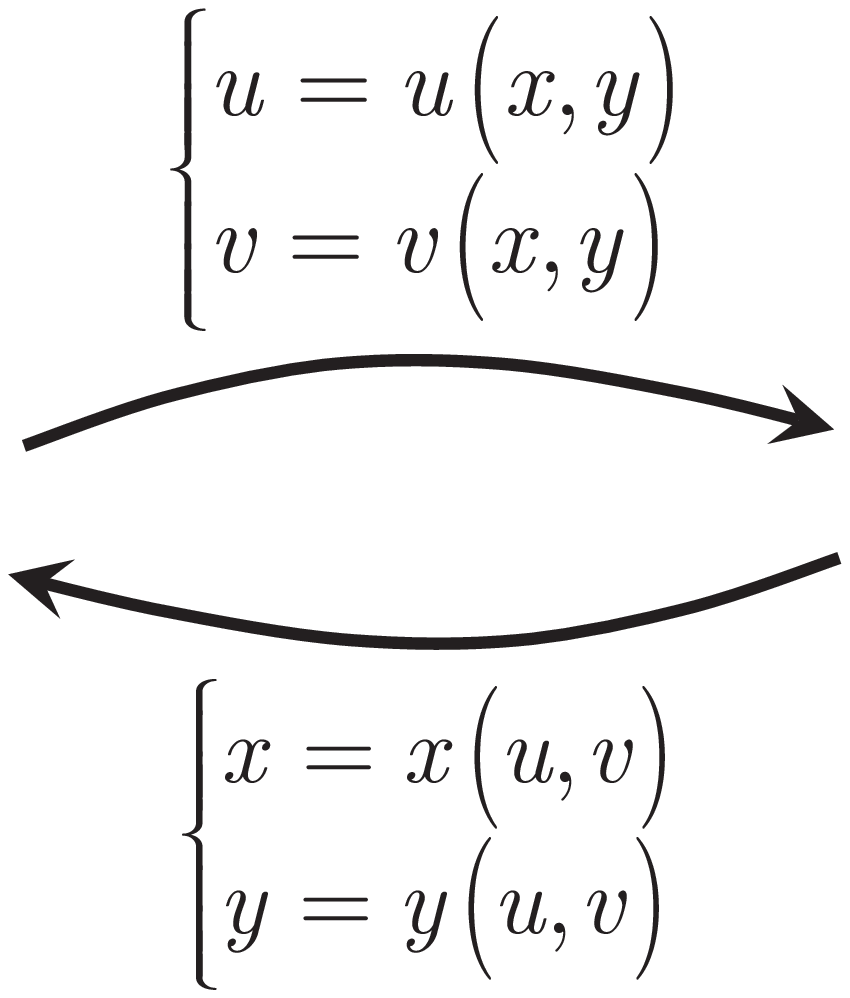}}
\scalebox{0.16}{\includegraphics{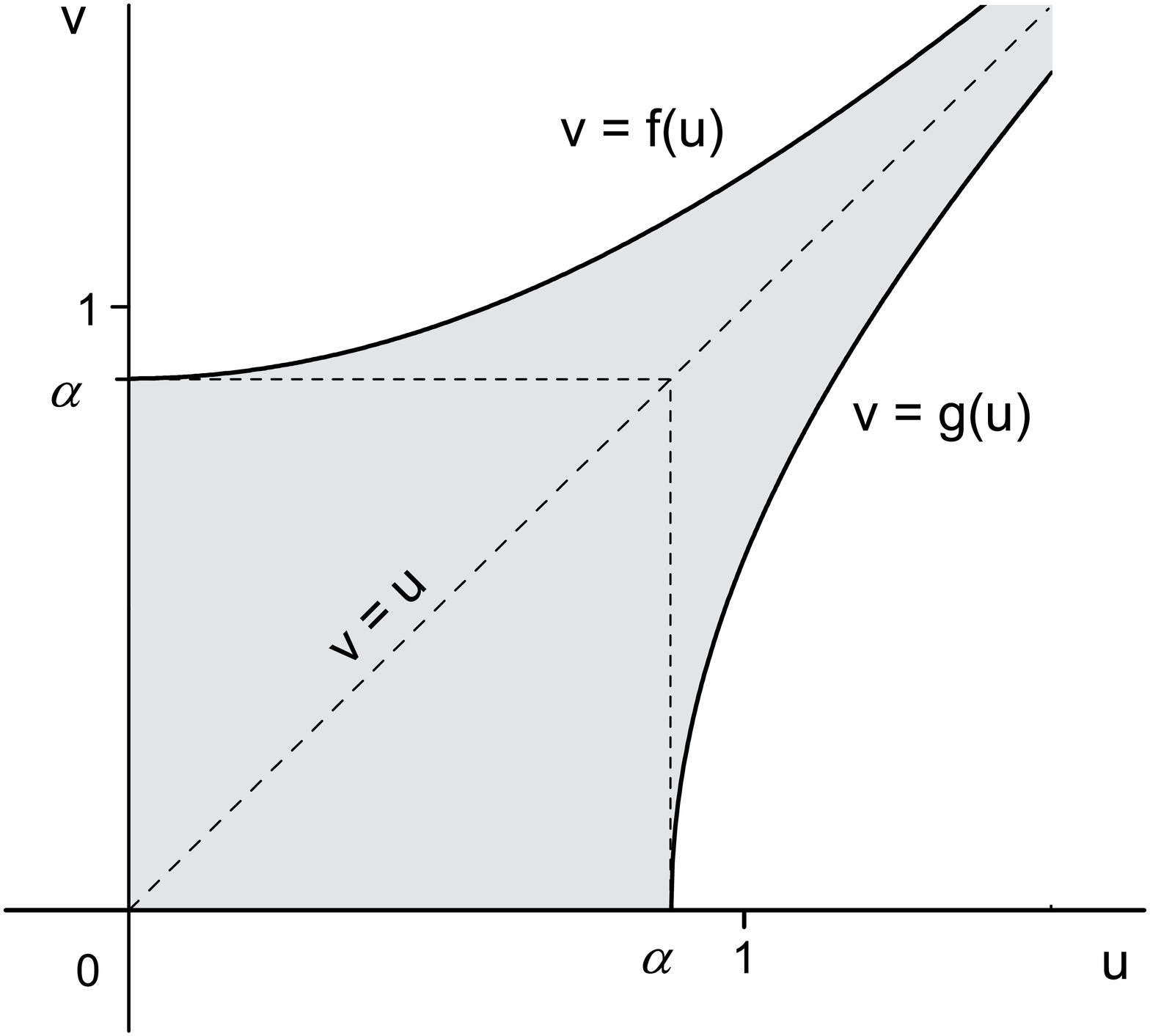}}
\caption{\label{fig:mapeamento}  Mapping the unit square in the $x y$-plane (left-hand side) onto the hyperbolic `quadrilateral' in the $u v$-plane (right-hand side). The symmetry about the diagonal $v=u$, represented by the inclined dashed line, is exploited in the text. Note that the hyperbolic region formed in the $u v$-plane is bounded by the coordinate axes and the curves $v=f(u)$ and $v=g(u)$. Both these functions, as well as the expressions for the change of variables indicated in the central part, are stated elsewhere in the text. According to Eq.~\eqref{eq:alfa}, $\alpha$ evaluates to $\: \sinh^{-1}(1) = \ln\!\left(1+\sqrt{2}\,\right) \approx 0.88137$.}
\end{center}
\end{figure}

Let us introduce the following \emph{hyperbolic version} of BCK change of variables for double integrals over the unit square:
\begin{equation}
x = \frac{\sinh{u}}{\cosh{v}} \quad \text{and} \quad y = \frac{\sinh{v}}{\cosh{u}} \, .
\label{eq:hyper}
\end{equation}
The corresponding Jacobi determinant is:
\begin{equation*}
J = \left|
  \begin{array}{cc}
    \dfrac{\cosh{u}}{\cosh{v}} & -\dfrac{\sinh{u} \, \sinh{v}}{\cosh^2{v}} \\
    \\
    -\dfrac{\sinh{u} \, \sinh{v}}{\cosh^2{u}} & \dfrac{\cosh{v}}{\cosh{u}} \\
  \end{array}
\right| = 1 -\tanh^2{u}\,\tanh^2{v} .
\label{eq:detJ}
\end{equation*}
As $0 \le x \, y \le 1$ in the unit square, then $1 - x^2\,y^2 \ge 0$, then $|J|$ is also equal to $1 -\tanh^2{u}\,\tanh^2{v}$. Therefore
\begin{eqnarray}
\int_0^1\!\!\!{\int_0^1{F(x,y) \, dx} \, dy} &=& {\int \! \! \! \int}_{\!\! S} \, F(x(u,v),y(u,v)) \: \left|J(u,v)\right| \: du \, dv \, \nonumber \\
&=& {\int \! \! \! \int}_{\!\! S} \, \tilde{F}(u,v) \: \left(1 -\tanh^2{u}\,\tanh^2{v}\right) \, du \, dv \, ,
\label{HypChange}
\end{eqnarray}
where $S$ is the cusped hyperbolic `quadrilateral' on the right-hand side of Fig.~\ref{fig:mapeamento}.  We are assuming that this change of variables is a $C^{1}$ diffeomorphism, which is reasonable since the functions $x=x(u,v)$ and $y=y(u,v)$ given in Eq.~\eqref{eq:hyper}, as well as
\begin{equation*}
u = u(x,y) = \cosh^{-1}{\sqrt{\frac{1+x^2}{1-x^2 \, y^2}}}
\end{equation*}
and
\begin{equation*}
v = v(x,y) = \cosh^{-1}{\sqrt{\frac{1+y^2}{1-x^2 \, y^2}}} \, ,
\end{equation*}
which define the \emph{inverse} change,\footnote{Note that the real function $\sinh^{-1}{x}$ is the inverse of the \emph{bijective} function $\sinh{x}$ for all $x \in \mathbb{R}$, as usual.  However, since the function $\cosh{x}$ is not bijective over the entire real domain, we can reduce its domain to the nonnegative reals and then \emph{define} its inverse as $\cosh^{-1}{x} := \ln{\left(x+\sqrt{x^2-1}\,\right)}$, $\forall \, x \ge 1$.}  are continuously differentiable in all points of their domains (see the hachured regions in Fig.~\ref{fig:mapeamento}), the only exception being point $(x,y) = (1,1)$.  As indicated there in Fig.~\ref{fig:mapeamento}, $S$ is bounded by the axes $u$ and $v$, the curve $v = f(u) = \sinh^{-1}{(\cosh{u})}$, above the diagonal $v=u$, and the curve $v = g(u) = \cosh^{-1}{(\sinh{u})}$, below the diagonal.  The parameter $\,\alpha\,$ can be readily evaluated by noting that it is the abscissa of the point where $v = g(u)$ intersects the $u$-axis. There in that point, one has $\,\cosh^{-1}{\!\left(\sinh{\alpha}\right)} = 0$, thus $\sinh{\alpha} = 1$. Since $\sinh^{-1}{x} = \ln{\left(x+\sqrt{x^2+1}\,\right)}$, $\forall \, x \in \mathbb{R}$, then
\begin{equation}
\alpha = \sinh^{-1}{(1)} = \ln\!{\left( 1+\sqrt{2}\, \right)} .
\label{eq:alfa}
\end{equation}

\section{New definite integrals and a dilogarithm identity}
Let us now apply our hyperbolic change of variables to the unit square integral $\,K$, as defined in Eq.~\eqref{eq:intK}.

\begin{teo}[First definite integral]
\label{teo1}
\begin{equation*}
\int_0^{\,\infty}\!{\left[\sinh^{-1}{\left(\cosh{u}\right)} -u\right] \, du} = \frac{\pi^2}{16} \, .
\end{equation*}
\end{teo}

\begin{prova}
From the BCK proof, one knows that $\int_0^1\!\!{\int_0^1{{\,1/(1-x^2 y^2)} \, dx} \, dy} = {\,\pi^2/8}$. By combining this result with~\eqref{HypChange}, one has
\begin{eqnarray*}
\frac{\pi^2}{8}  &=& {\int \! \! \! \int}_{\!\! S} \, \frac{1}{1 -\tanh^2{u}\,\tanh^2{v}} \, \left(1 -\tanh^2{u}\,\tanh^2{v}\right) \: du \, dv \nonumber \\
&=& {\int \! \! \! \int}_{\!\! S} \: 1 \: du \, dv \, = \, \text{area of~} S .
\end{eqnarray*}
The area of the domain $S$ can be equalized to a single definite integral by noting that $S$ is symmetric with respect to the diagonal $v=u$ (see Fig.~\ref{fig:mapeamento}). Then, by taking only the upper half of $S$ into account, one has
\begin{equation*}
\text{area of~} S = 2 \int_0^{\,\infty}{\!\!\!{\int_u^{f(u)} \: 1 \: dv} \, du} = 2 \int_0^{\,\infty}\!{\left[\sinh^{-1}{\left(\cosh{u}\right)} -u\right] \, du} \, ,
\end{equation*}
which leads to the desired result. \qquad \qquad \qquad \qquad \qquad \qquad \qquad \qquad \qquad \quad $\Box$
\end{prova}

\begin{teo}[Second definite integral]
\label{teo2}
\begin{equation*}
\int_\alpha^{\,\infty}\!{\left[u -\cosh^{-1}{\left(\sinh{u}\right)}\right] \, du} = \frac{\pi^2}{16} -\frac12\,\ln^2\!{\left(1+\sqrt{2}\,\right)} \, .
\end{equation*}
\end{teo}

\begin{prova}
Similarly to the previous proof, by taking the lower half of $S$ into account and noting that it is bounded by the $u$-axis, the diagonal $v=u$, and the curve $v = g(u)$, as seen in Fig.~\ref{fig:mapeamento}, one has, from Eq.~\eqref{HypChange}, that:
\begin{equation*}
\text{area of~} S = 2 \, \left[\frac{\alpha^2}{2} + \int_\alpha^{\,\infty}{\!\!\!{\int_{g(u)}^u \: 1 \: dv} \, du} \right] = \alpha^2 + 2 \, \int_\alpha^{\,\infty}\!{\left[u - \cosh^{-1}{\left(\sinh{u}\right)} \right] \, du} \, .
\end{equation*}
Again, since the area of $S$ is ${\,\pi^2/8}$, one finds
\begin{equation*}
\int_\alpha^{\,\infty}\!{\left[u - \cosh^{-1}{\left(\sinh{u}\right)} \right] \, du} = \frac12 \left(\frac{\pi^2}{8} -\alpha^2\right) = \frac{\pi^2}{16} -\frac{\alpha^2}{2} \, .
\end{equation*}
\begin{flushright} $\Box$ \end{flushright}
\end{prova}

\begin{teo}[Third definite integral]
\label{teo3}
\begin{equation*}
\int_{\,\alpha/2}^{\,\infty}{\,\ln{\left(\tanh{z}\right)} \: dz} = \frac14 \, \ln^2{\!\left(1+\sqrt{2}\,\right)} -\frac{\pi^2}{16} \, .
\end{equation*}
\end{teo}

\begin{prova}
From Theorem~\ref{teo1}, we know that the area of $S$ evaluates to
\begin{equation*}
2 \int_0^{\,\infty}{\!\!{\int_u^{\,f(u)} 1 \: dv} \: du} = \frac{\pi^2}{8} \, ,
\end{equation*}
where $f(u)=\sinh^{-1}{(\cosh{u})}$. Let us rotate the coordinate axes by the angle $\pi/4$~rad, as in Apostol's proof. Since ``pure'' rotations about the origin always yield an unitary Jacobi determinant, one has:
\begin{eqnarray}
\frac{\pi^2}{8} &=& 2 \int_0^{\,\alpha/\sqrt{2}}{\!\!\int_0^X{ 1 \: dY} \, dX} +2 \int_{\alpha/\sqrt{2}}^\infty{\,\int_0^{\,h(X)}{1 \; dY}\,dX} \nonumber \\
&=& 2 \int_0^{\,\alpha/\sqrt{2}}{X \, dX} +2 \int_{\alpha/\sqrt{2}}^\infty{\,h(X) \, dX} \nonumber \\
&=& \frac{\alpha^2}{2} + 2 \int_{\alpha/\sqrt{2}}^\infty{\,h(X) \: dX} ,
\label{eq:aux}
\end{eqnarray}
where
\begin{equation*}
u = \frac{X-Y}{\sqrt{2}} \; , \quad v = \frac{X+Y}{\sqrt{2}} \, ,
\end{equation*}
and $h(X)$ is determined by changing the variables $u$ and $v$ on the equation $v=f(u)$ to the new variables $X$ and $Y$, as follows. Since $v = f(u) = \sinh^{-1}(\cosh{u})$, then $\sinh{v} = \cosh{u}$, so
\begin{equation*}
\sinh{\left( \frac{X+Y}{\sqrt{2}} \right)} = \cosh{\left( \frac{X-Y}{\sqrt{2}} \right)} \, .
\end{equation*}
By making use of the addition and subtraction formulas for hyperbolic functions, as well as of the substitutions $X' := X/\sqrt{2}$ and $Y' := Y/\sqrt{2}$, the above equality expands to
\begin{equation*}
\sinh{X'} \: \cosh{Y'} + \sinh{Y'} \: \cosh{X'} = 
\cosh{X'} \: \cosh{Y'} - \sinh{X'}  \: \sinh{Y'} \, .
\end{equation*}
By dividing both sides by $\cosh{X'} \: \cosh{Y'}$, one finds
\begin{equation*}
\tanh{X'} + \tanh{Y'} = 1 -\tanh{X'} \cdot \tanh{Y'} \, ,
\end{equation*}
which can be written as
\begin{equation}
\tanh{Y'} = \frac{1 -\tanh{X'} }{1 +\tanh{X'} } \, .
\label{eq:dem1}
\end{equation}
Now, for simplicity, let us define $\, t := \tanh{X'}$. After some algebra, Eq.~\eqref{eq:dem1} becomes
\begin{equation*}
Y' = \frac12 \: \ln{\left( \frac{1 + \dfrac{1-t}{1+t}}{1 -\dfrac{1-t}{1+t}} \right)} = \frac12 \: \ln{\left(\frac{1}{t}\right)} = -\frac12 \: \ln{t} \, ,
\end{equation*}
which means that
\begin{equation*}
Y = h(X) = -\frac{\sqrt{2}}{2} \; \ln{\!\left[ \tanh{\!\left(\frac{X}{\sqrt{2}}\right)} \right]} .
\end{equation*}
By putting this function on the last integral of Eq.~\eqref{eq:aux}, one has
\begin{equation*}
\frac{\pi^2}{8} = \frac{\alpha^2}{2} - \sqrt{2} \, \int_{\alpha/\sqrt{2}}^\infty{\:\ln{\!\left[\tanh{\left(\frac{X}{\sqrt{2}}\right)}\right]} \, dX} .
\label{eq:aux2}
\end{equation*}
The simple substitution $\,z = {\,X/\sqrt{2}}\,$ promptly reduces this equation to
\begin{equation}
\int_{\alpha/2}^{\,\infty}{\ln{\left(\tanh{z}\right)} \: dz} = \, \frac{\alpha^2}{4} -\frac{\pi^2}{16} \, .
\label{eq:aux3}
\end{equation}
\begin{flushright} $\Box$ \end{flushright}
\end{prova}

After this proof, we make the substitution $\,t = \tanh{z}\,$ in the integral in Eq.~\eqref{eq:aux3}, above. When the result is expressed in terms of the dilogarithm function, an interesting identity arises.

\begin{teo}[Dilogarithm identity]
\label{teo4}
Let us define $\,\mathrm{Li}_2{(z)} := \sum_{k=1}^\infty{{\,z^k/k^2}}$, which is valid for all complex $z$ with $|z| \le 1$~\emph{\cite{Lewin}}. ~Then
\begin{equation*}
\mathrm{Li}_2{\left(\sqrt{2}-1\right)} +\mathrm{Li}_2{\left( 1 -\frac{1}{\sqrt{2}} \right)} = \frac{\pi^2}{8} -\frac{\ln^2{\!\left(1+\sqrt{2}\,\right)}}{2} -\frac18 \, \ln^2{2} \, .
\label{eq:fim}
\end{equation*}
\end{teo}

\begin{prova}
By applying the substitution $\,t = \tanh{z}\,$ ($dt = \mathrm{sech}^2{z}\,dz$) to the definite integral of Theorem~\ref{teo3}, in the previous section, one finds
\begin{equation}
\frac{\alpha^2}{4} -\frac{\pi^2}{16} = \int_{\tanh{\!({\alpha/2})}}^1{\!\frac{\ln{t}}{1-t^2} \; dt} 
= \int_{\sqrt{2}-1}^1{\frac{\ln{t}}{1-t^2} \; dt} \, ,
\label{eq:ln}
\end{equation}
since
\begin{equation*}
\tanh{\!\left(\frac{\alpha}{2}\right)} = \tanh{\!\left(\ln{\sqrt{1+\sqrt{2}}}\right)} = \frac{\sqrt{1+\sqrt{2}} -\dfrac{1}{\sqrt{1+\sqrt{2}}}}{\sqrt{1+\sqrt{2}} +\dfrac{1}{\sqrt{1+\sqrt{2}}}} = \frac{\sqrt{2}}{2+\sqrt{2}} = \sqrt{2}-1 .
\end{equation*}
By expanding the integrand in Eq.~\eqref{eq:ln} in partial fractions, one has
\begin{equation}
\frac{\alpha^2}{4} -\frac{\pi^2}{16}
= \frac12 \, \int_{\sqrt{2}-1}^1{\frac{\ln{t}}{1+t} \: dt} -\frac12 \, \int_1^{\sqrt{2}-1}{\!\frac{\ln{t}}{1-t} \: dt} \, .
\label{eq:2lns}
\end{equation}
These integrals can be reduced to special values of the dilogarithm function if one adopts the more general integral definition for this function~\cite{Lewin}, namely
\begin{equation*}
\mathrm{Li}_2{(z)} := -\int_0^{z}{\frac{\ln{(1-s)}}{s} \; ds} \, ,
\end{equation*}
valid for $z \in \mathbb{C} \backslash [1,\infty)$.\footnote{We are considering here the principal branch of the dilogarithm function $\mathrm{Li}_2{(z)}$, defined by taking the principal branch of $\ln{z}$, for which it has a cut along the negative real axis, with $|\arg{(z)}| < \pi$. This defines the principal branch of $\mathrm{Li}_2{(z)}$ as a single-valued function in the complex plane cut along the real axis, from $1$ to $+\infty$.}  Note that, by expanding the logarithm in powers of $z$ in this integral definition, one finds $\mathrm{Li}_2{(z)} = \sum_{k=1}^\infty{{\,z^k/k^2}}$, which agrees with the basic definition stated in the theorem for all complex $z$ with $|z| \le 1$~\cite{Lewin}. By taking into account this summation form, it is clear that $\mathrm{Li}_2{(1)} = {\,\pi^2/6}$, which implies that
\begin{equation*}
\mathrm{Li}_2{(z)}  = \frac{\pi^2}{6} -\,\int_1^{\,z}{\frac{\ln{(1-s)}}{s} \: ds} \, .
\label{eq:defLi2}
\end{equation*}
From this integral form, it is easy to deduce that
\begin{equation*}
\int_1^{\,z}{\frac{\ln{t}}{1-t} \: dt} = \mathrm{Li}_2{(1-z)}
\end{equation*}
and
\begin{equation*}
\int_z^{\,1}{\frac{\ln{t}}{1+t} \: dt} = -\mathrm{Li}_2{(-z)} -\ln{z}\,\ln{(z+1)} -\frac{\pi^2}{12} \, .
\end{equation*}

By putting $z=\sqrt{2}-1$ in these integrals and then substituting the resulting dilogarithms in Eq.~\eqref{eq:2lns}, one finds
\begin{eqnarray}
\frac{\alpha^2}{4} -\frac{\pi^2}{16} &=& \frac12 \, \left[ -\mathrm{Li}_2{\left(1-\sqrt{2}\right)} -\ln{\!\left(\sqrt{2}-1\right)}\,\ln{\sqrt{2}} -\frac{\pi^2}{12} \right] -\frac12 \, \mathrm{Li}_2{\left(2-\sqrt{2}\right)} \nonumber \\
&=& -\frac{\mathrm{Li}_2{\left(1-\sqrt{2}\,\right)}}{2} -\frac{\mathrm{Li}_2{\left(2-\sqrt{2}\,\right)}}{2} -\frac14\,\ln{\!\left(\sqrt{2}-1\right)}\,\ln{2} -\frac{\pi^2}{24} \, .
\label{eq:auxLi2}
\end{eqnarray}
Now, by applying the Euler reflection formula for $\mathrm{Li}_2{(z)}$ about $z=\frac12$, which reads (see Eq.~(1.11) in Ref.~\cite{Lewin})
\begin{equation*}
\mathrm{Li}_2{(1-z)} = -\mathrm{Li}_2{(z)} -\ln{z}\,\ln{(1-z)} +\frac{\pi^2}{6} \, ,
\end{equation*}
to $\mathrm{Li}_2{\left(2-\sqrt{2}\,\right)}$ in Eq.~\eqref{eq:auxLi2}, one finds
\begin{equation}
\mathrm{Li}_2{\left(\sqrt{2}-1\right)} - \mathrm{Li}_2{\left(1 -\sqrt{2}\,\right)} = \frac{\alpha^2}{2} -\frac{\pi^2}{8} + \frac{\pi^2}{4} -\ln^2{\!\left(\sqrt{2}-1\right)} .
\label{eq:tmp}
\end{equation}
By applying the Landen's formula, namely (see Eq.~(1.12) in Ref.~\cite{Lewin})
\begin{equation*}
\mathrm{Li}_2{(z)} = -\mathrm{Li}_2{\left( \frac{z}{z-1} \right)} -\frac12 \, \ln^2{(1-z)} \, ,
\end{equation*}
which is valid for all complex $z$, $z \not \in [1,\infty)$, to $\,\mathrm{Li}_2{\left(1 -\sqrt{2}\,\right)}$, one has
\begin{equation*}
\mathrm{Li}_2{\left(1 -\sqrt{2}\,\right)} = -\mathrm{Li}_2{\left( 1 -\frac{1}{\sqrt{2}} \right)} -\frac18 \, \ln^2{2} \, .
\end{equation*}
Finally, by substituting this result on Eq.~\eqref{eq:tmp}, one finds
\begin{equation*}
\mathrm{Li}_2{\left(\sqrt{2}-1\right)} +\mathrm{Li}_2{\left( 1 -\frac{1}{\sqrt{2}} \right)} = \frac{\pi^2}{8} -\frac{\alpha^2}{2} -\frac18 \, \ln^2{2} \, .
\label{eq:fim2}
\end{equation*}
\begin{flushright} $\Box$ \end{flushright}
\end{prova}

It should be mentioned that neither Mathematica (release 7) nor Maple (release 13) are able to generate any closed-form expression for this sum.  Indeed, this identity is not found in standard bibliographic sources for dilogarithms~\cite{Lewin,Zagier}, though it resembles the identity
\begin{equation*}
\mathcal{L}{\left(\sqrt{2}-1\right)} +\mathcal{L}{\left( 1 -\frac{1}{\sqrt{2}} \right)} = \frac34 \: ,
\label{eq:rogers}
\end{equation*}
where $\:\mathcal{L}{(x)} := \dfrac{6}{\pi^2} \left[ \mathrm{Li}_2{(x)} +\frac12 \, \ln{x} \: \ln{(1-x)} \right]\,$ is the \emph{normalized Roger's dilogarithm}, as found by Bytsko (1999)~\cite{Bytsko}.

It seems plausible to extend the results obtained here in this paper to other known unit square integrals, such as $\,\int_0^1\!\!{\int_0^1{\ln{(x y)}/(1-x^2 y^2) \, dx} \, dy} = -\frac74 \, \zeta{(3)}\,$, $\,\int_0^1\!\!{\int_0^1{\ln{x}/(1-x^2 y^2) \, dx} \, dy} = -\frac78 \, \zeta{(3)}\,$, $\,\int_0^1\!\!{\int_0^1{\ln{(x y)}/(1-x y) \, dx} \, dy} = -2 \, \zeta{(3)}\,$, $\,\int_0^1\!\!{\int_0^1{\ln{(1-x y)}/(1-x y) \, dx} \, dy} = -\zeta{(3)}\,$, and many others (see, e.g., Ref.~\cite{Wolfram}). Similarly, it should also be interesting to include a third dimension in the hyperbolic change of variables in order to investigate some triple integrals over the unit cube, such as
\begin{equation*}
\int_0^1\!\!{\int_0^1\!\!{\int_0^1{\frac{1}{1 -x y z} \: dx} \, dy} \, dz} = \zeta{(3)} \quad \text{and} \; \int_0^1\!\!{\int_0^1\!\!{\int_0^1{\frac{1}{1 -x^2 y^2 z^2} \: dx} \, dy} \, dz} = \frac78 \, \zeta{(3)} \, .
\end{equation*}
Extensive work in both directions is already ongoing by the author.

\section*{Acknowledgments}
The author acknowledges a postdoctoral fellowship from CNPq (Brazilian agency) during the course of this work.

\newpage

\section*{Figures}

\end{document}